\documentstyle[12pt]{article}
\begin{document}
\newpage
\pagestyle{empty}
\setcounter{page}{0}
\renewcommand{\thesection}{\arabic{section}}
\renewcommand{\theequation}{\thesection.\arabic{equation}}
\newcommand{\sect}[1]{\setcounter{equation}{0}\section{#1}}
\newfont{\twelvemsb}{msbm10 scaled\magstep1}
\newfont{\eightmsb}{msbm8}
\newfont{\sixmsb}{msbm6}
\newfam\msbfam
\textfont\msbfam=\twelvemsb
\scriptfont\msbfam=\eightmsb
\scriptscriptfont\msbfam=\sixmsb
\catcode`\@=11
\def\Bbb{\ifmmode\let\next\Bbb@\else
  \def\next{\errmessage{Use \string\Bbb\space only in math mode}}\fi\next}
\def\Bbb@#1{{\Bbb@@{#1}}}
\def\Bbb@@#1{\fam\msbfam#1}
\newfont{\twelvegoth}{eufm10 scaled\magstep1}
\newfont{\tengoth}{eufm10}
\newfont{\eightgoth}{eufm8}
\newfont{\sixgoth}{eufm6}
\newfam\gothfam
\textfont\gothfam=\twelvegoth
\scriptfont\gothfam=\eightgoth
\scriptscriptfont\gothfam=\sixgoth
\def\frak{\frak@}
\def\frak@#1{{\fam\gothfam{{#1}}}}
\def\frak@@#1{\fam\gothfam#1}
\catcode`@=12
%
%
%
\def\CC{{\Bbb C}}
\def\NN{{\Bbb N}}
\def\QQ{{\Bbb Q}}
\def\RR{{\Bbb R}}
\def\ZZ{{\Bbb Z}}
\def\cA{{\cal A}}          \def\cB{{\cal B}}          \def\cC{{\cal C}}
\def\cD{{\cal D}}          \def\cE{{\cal E}}          \def\cF{{\cal F}}
\def\cG{{\cal G}}          \def\cH{{\cal H}}          \def\cI{{\cal I}}
\def\cJ{{\cal J}}          \def\cK{{\cal K}}          \def\cL{{\cal L}} 
\def\cM{{\cal M}}          \def\cN{{\cal N}}          \def\cO{{\cal O}}
\def\cP{{\cal P}}          \def\cQ{{\cal Q}}          \def\cR{{\cal R}} 
\def\cS{{\cal S}}          \def\cT{{\cal T}}          \def\cU{{\cal U}}
\def\cV{{\cal V}}          \def\cW{{\cal W}}          \def\cX{{\cal X}}
\def\cY{{\cal Y}}          \def\cZ{{\cal Z}}
\def\qed{\hfill \rule{5pt}{5pt}}
\def\id{\mbox{id}}
\def\ggo{{\frak g}_{\bar 0}}
\def\uqggo{\cU_q({\frak g}_{\bar 0})}
\def\uqggp{\cU_q({\frak g}_+)}
\def\typeA{{\em type $\cA$}}
\def\typeB{{\em type $\cB$}}
\newtheorem{lemma}{Lemma}
\newtheorem{prop}{Proposition}
\newtheorem{theo}{Theorem}
\newtheorem{Defi}{Definition}
$$
\;
$$
$$
\;
$$
$$
\;
$$
$$
\;
$$
\vfill
\vfill
\begin{center}

{\LARGE {\bf {\sf 
Reflection Equation Algebra of a $(h,w)$-deformed Oscillator
}}} \\[2cm]

{\large Boucif Abdesselam\footnote{boucif@orphee.polytechnique.fr} and 
Ranabir Chakrabarti\footnote{Permanent address: 
Department of Theoretical Physics, University of Madras, Guindy Campus, 
Madras-600025, India}}

\smallskip 

\smallskip 

\smallskip 

{\em  \footnote{Laboratoire Propre du CNRS UPR A.0014}Centre de Physique 
Th\'eorique, Ecole Polytechnique, \\
91128 Palaiseau Cedex, France.}

\end{center}

\smallskip 

\smallskip

\vfill

\begin{abstract}
 We consider the reflection equation algebra for a finite dimensional 
$R$-matrix for the $(h,w)$-deformed Heisenberg algebra ${\cal U}_{h,w}(h(4))$.
A representation of the reflection matrix $K$ is constructed using the matrix 
generators $L^{(\pm)}$ of the ${\cal U}_{h,w}(h(4))$ algebra. A series of 
representations of the $K$-matrix then may be generated by using the 
coproduct rules of the ${\cal U}_{h,w}(h(4))$ algebra. The complementary 
condition necessary for combining two distinct solutions of the reflection 
equation algebra yields the braiding relations between these two sets of 
generators. This may be thought as a generalization of Bose-Fermi statistics 
to braiding statistics, which them may be used to provide a new braided 
colagebraic 
structure to a Hopf algebra generated by the elements of the matrix $K$.
The reflection equation algebra and the braided exchange properties are found 
to depend on both deformation parameters $h$ and $w$. 
\end{abstract}

\vfill
\vfill

\rightline{CPTH-S526.0697}
\rightline{q-alg/9706032}
\rightline{June 97}

$$
\;
$$
$$
\;
$$

\newpage 

\pagestyle{plain}

\sect{Introduction}

The representation theory of quantum algebras has led to various deformed 
oscillator algebras [1-7]. These studies may lead to field theories, where 
excitations obeying braiding statistics [4,8] may be discussed. Other 
applications of deformed Heisenberg algebra include the description [9] 
of a class of exactly solvable potentials in terms of dynamical symmetries, 
the theory of link invariants [10] and the theory of $q$-special functions 
[11]. At an appropriate contraction limit of the quantum algebra 
${\cal U}_q(sl(2))$, a `standard' deformed Heisenberg algebra 
${\cal U}_h(h(4))$ and its universal ${\cal R}$-matrix has been obtained 
[12]. A different `nonstandard' quantization of the Heisenberg algebra 
${\cal U}_w(h(4))$ has also been obtained [13] as a contraction limit of the 
Jordanian deformation of $sl(2)$ algebra [14]. A distinction between these 
two deformations appears in their properties 
under a rescaling of the Heisenberg generators: $A_{\pm} \rightarrow 
\lambda^{\pm 1}A_{\pm}$. While the `standard' ${\cal U}_h(h(4))$ 
algebra remains invariant under this rescaling, the `nonstandard' algebra 
${\cal U}_w(h(4))$ does so only if we also rescale the deformation 
parameter $w\rightarrow \lambda^{-1}w$. In this sense, the algebra 
${\cal U}_w(h(4))$ has close kinship with $\kappa$-Poincar\'e 
algebra [15] and deformed conformal algebra [16], where the introduction of 
mass-like deformation parameter leads to the appearance of the fundamental 
mass on basic geometrical level. Combining the features of the `standard' 
and `non-standard' deformation procedures, a new two-parametric quantization 
of the Heisenberg algebra ${\cal U}_{h,w}(h(4))$ has recently been 
obtained [17]. The full Hopf structure of the ${\cal U}_{h,w}(h(4))$ algebra, 
its universal ${\cal R}$-matrix and the corresponding matrix pseudo-group 
have also been discussed. From the point of view of applicability in concrete
physical models, the algebra ${\cal U}_{h,w}(h(4))$ with multiparameter 
$(h,w)$-deformation may be of interest.

\smallskip

In another development, the reflection equation (RE) algebra connected to 
the quantum group has attracted wide interest. The RE was introduced in [18] 
to describe factorized scattering on a half-line, and the related algebra 
soon found quite different applications in quantum current algebras [19],  
integrable models with non-periodic boundary conditions [20] and in the 
description [21] of braid group representations on a handlebody. A 
generalization of the RE algebra for the scattering of lines moving in a half 
plane touching the boundary in $2+1$-dimension has been made [22]. Recently, 
there has been renewed interest in this topic following investigations in 
condensed matter physics in which boundaries play a significant role. The 
boundary states, which correspond to the lowest lying energy states in 
integrable field 
theories and statistical mechanics models with boundaries have been explored
[23,24]. 
In an attempt to establish a second quantized approach to this problem, an 
oscillator algebra and the associated Fock space satisfying reflection 
boundary condition and obeying generalized statistics has been studied [25]. 
In this connection it is of interest to study the RE algebra corresponding 
to the $R$-matrix of the doubly deformed Heisenberg algebra 
${\cal U}_{h,w}(h(4))$. The braided group approach [26,27] to the RE 
algebras takes the point of view that in the exchanges between two copies
of the RE algebras, 
the usual transposition map is replaced by a more general braiding. From this 
point of view, the braided group corresponding to the singly deformed 
Heisenberg algebra ${\cal U}_h(h(4))$ with the `standard' deformation 
parameter $h$ has been studied [28]. In the 
present context, it may be useful to investigate the braided exchange
properties described by the RE algebra of the $(h,w)$-deformed oscillator. 
In particular, the role of the `nonstandard'  
deformation parameter $w$ in the structure of the braiding statistics is worth
exploring. This may serve as a pointer towards understanding the role of the 
dimensional mass-like deformation parameters in the context of RE algebras 
corresponding to the $\kappa$-deformed Poincar\'e algebra [15] and deformed 
conformal algebras [16].
  
\smallskip

In this article, we study the spectral parameter independent form of the RE 
algebra related to the quantum algebra ${\cal U}_{h,w}(h(4))$. In Section 2, 
we review the Hopf structure of the doubly deformed algebra 
${\cal U}_{h,w}(h(4))$ and recast it using the FRT prescription [29]. The 
extended RE algebra related to its $R$-matrix is described in Section 3. 
The braiding properties of the two copies of the generators of the RE 
algebra depend on both the deformation parameters $(h,w)$. We conclude 
in Section 4.

\sect{The Hopf Algebra ${\cal U}_{h,w}(h(4))$ and its 
FRT Construction} 

We start by enlisting the full Hopf structure [17] of the universal enveloping 
algebra ${\cal U}_{h,w}(h(4))$ generated by the elements $A_{\pm}$,
$N$ and $E$. The algebra reads  
\begin{eqnarray}
  &&  [N,A_{+}]= {e^{wA_{+}}-1\over w},\qquad\;\;\qquad\qquad  
[N,A_-]= - A_-, \nonumber \\
 && [A_{-},A_{+}]={\sinh (hE)\over h}e^{w A_{+}}, \qquad\qquad [E,\;\bullet\;]=0,
\qquad
\end{eqnarray}
where $E$ is a central element of the algebra. Another nonlinear central 
element $C$ exists [17] and may be regarded as the Casimir element of the 
algebra (2.1):
\begin{eqnarray}
C=  {\sinh{ h E} \over h}N + {e^{-w A_{+}} -1 \over 2\;w} A_{-}+ 
A_{-} {e^{-w A_{+}} -1 \over 2\;w}.  
\end{eqnarray}
The coalgebraic structure of ${\cal U}_{h,w}(h(4))$ is given by [17]
\begin{eqnarray}
&& \bigtriangleup(A_{+})=A_{+} \otimes 1 + 1 \otimes A_{+}, \nonumber \\
&& \bigtriangleup(A_{-})=A_{-} \otimes e^{hE}e^{w A_{+}}+e^{-hE}\otimes A_{-}
+we^{-hE}\;N\otimes {\sinh (hE) \over h}e^{w A_{+}}, \nonumber \\
&& \bigtriangleup(N)=N\otimes  e^{w A_{+}}+ 1\otimes N, \qquad \qquad 
\bigtriangleup(E)=E\otimes 1 + 1\otimes E,\\
&& \varepsilon(A_{+})=\varepsilon(A_{-})=\varepsilon (N)= \varepsilon(E)= 0,
\qquad\qquad\qquad \\
&& S(A_{+})=-A_{+},\qquad\qquad S(A_{-})=-A_{-}e^{-w A_{+}}+w\;{\sinh (hE) 
\over h}N\;e^{-w A_{+}}, \nonumber \\
&& S(N)=-Ne^{-w A_{+}}, \qquad\qquad S(E)=-E.
\end{eqnarray}
Here we note that there exists an invertible map of the algebra (2.1) on the 
undeformed Heisenberg algebra:
\begin{eqnarray}
&& a_{+}={1-e^{-w A_{+}}\over w}, 
\qquad a_{-}= {hE\over \sinh(hE)}A_{-},\qquad  n= N, \qquad  e= E,  
\end{eqnarray}
where the generators $(a_{\pm}, n,e)$ are classical in nature: 
\begin{eqnarray}   
&& [n,a_{\pm}]=\pm a_{\pm},\qquad\qquad [a_{-},a_{+}]=e,
 \qquad\qquad [e,\;\bullet\;]=0. \qquad\qquad 
\end{eqnarray}
The universal ${\cal R}$-matrix of ${\cal U}_{h,w}(h(4))$ 
algebra reads [17]
\begin{eqnarray} 
{\cal R}= e^{-w A_{+}\otimes N}
e^{-2h\;E\otimes N}\;\exp \biggl(2h\;e^{hE}A_{-}
\otimes \biggl({1-e^{-w A_{+}}\over w}\biggr)\biggr)
e^{w N \otimes A_{+}}.
\end{eqnarray}

Following the FRT prescription [29], we now recast the Hopf structure of the 
${\cal U}_{h,w}(h(4))$ algebra in order to explicitly obtain the 
Lax operators $L^{(\pm)}$, which may be used to construct a hierarchy 
of solutions [21] of the RE algebra. To this end, we first obtain the conjugate
universal ${\tilde {\cal R}}$ matrix that satisfies Yang-Baxter equation and 
is defined as ${\tilde {\cal R}}={\cal R}_{+}^{-1}$, where ${\cal R}_{+}=
P{\cal R}P$ and $P$ is the transposition operator:
\begin{eqnarray} 
{\tilde {\cal R}}=e^{-w A_{+}\otimes N} \;\exp \biggl(-2h\;
\biggl({1-e^{-w A_{+}}\over w}\biggr)\otimes e^{hE}A_{-}\biggr)
e^{2h\;N\otimes E} e^{w N \otimes A_{+}}.
\end{eqnarray}
  
A real $3\times 3$-matrix representation of the algebra (2.1) remains 
undeformed and reads as follows:
\begin{eqnarray}
&&\pi_{3}(A_{+})=\pmatrix{0 & 0& 0\cr
 0 & 0& 1\cr 0 & 0& 0\cr}, \qquad\qquad \pi_{3}(A_{-})=\pmatrix{0 & 1& 0\cr
 0 & 0& 0\cr 0 & 0& 0\cr},  \nonumber \\
&& \pi_{3}(N)=\pmatrix{0 & 0& 0\cr
 0 & 1& 0\cr 0 & 0& 0\cr}, \qquad\;\;\qquad \pi_{3}(E)=\pmatrix{0 & 0& 1\cr
 0 & 0& 0\cr 0 & 0& 0\cr}.
\end{eqnarray}
Using the representation (2.10) in one of the sectors of the tensor products 
in the expressions of ${\cal R}$ and ${\tilde {\cal R}}$ in (2.8) and (2.9) 
respectively, the matrix operators $L^{(\pm)}$ may be directly read: 
\begin{eqnarray}
&& \displaystyle (\pi_{3}\otimes \id){\cal R}=
\pmatrix{1&&\displaystyle {2h\over w}(e^{wA_{+}}-1)& &-2hN\cr
 && && \cr
 0 && \displaystyle  e^{wA_{+}} &&\displaystyle   -wN \cr 
 && && \cr
0 && 0&& 1\cr}=S(L^{(-)}), \nonumber \\
&& (\id \otimes \pi_{3}){\cal R}=
\pmatrix{1&&  0 & & 0\cr
 && && \cr
 0 && \displaystyle  e^{-2hE-wA_{+}} &&
\displaystyle   2he^{-hE-wA_{+}}A_{-}+we^{-2hE-wA_{+}}N\cr 
 && && \cr
0 && 0&& 1\cr}=L^{(+)}, \nonumber \\
&& (\pi_{3}\otimes \id){\tilde {\cal R}}=
\pmatrix{1&& 0 & & 0\cr
 && && \cr
 0 && \displaystyle  e^{2hE+wA_{+}} && \displaystyle  -2he^{hE}A_{-}-wN \cr 
 && && \cr
0 && 0&& 1\cr}=S(L^{(+)}) \nonumber \\
&& (\id\otimes\pi_{3} ){\tilde {\cal R}}=
\pmatrix{1&&\displaystyle  -{2h\over w}(1-e^{-wA_{+}}) & &
\displaystyle  2he^{-wA_{+}}N\cr
 && && \cr
 0 && \displaystyle  e^{-wA_{+}} && \displaystyle  we^{-wA_{+}}N \cr 
 && && \cr
0 && 0&& 1\cr}=L^{(-)}.
\end{eqnarray}
The ${\cal R}$-matrix in the representation (2.10), defined as 
$R=(\pi_{3} \otimes\pi_{3} ){\cal R}$ is given by
\begin{eqnarray}
R=1\otimes 1+2h\;(e_{12}\otimes e_{23}-e_{13}\otimes e_{22})+w
\; (e_{22}\otimes e_{23}-e_{23}\otimes e_{22})
\end{eqnarray}
where $(e_{ij})_{kl}=\delta_{ik}\;\delta_{jl}$ and the indices $(i,j,k,l)\in 
(1,2,3)$. The transposed matrix $R_{+}=PRP$ and the conjugate
matrix ${\tilde R}=R_{+}^{-1}$ may be obtained readily. 
The algebra (2.1) may now be reformulated in terms of the matrix generators 
$L^{(\pm)}$: 
\begin{eqnarray}
R_{+}L_{1}^{(\epsilon_1)}L_{2}^{(\epsilon_2)}=L_{2}^{(\epsilon_2)}
L_{1}^{(\epsilon_1)}R_{+}
\end{eqnarray}
where $L_{1}^{(\epsilon)}=L^{(\epsilon)} \otimes 1$, $L_{2}^{(\epsilon)}
=1\otimes L^{(\epsilon)}$ and $(\epsilon_1,\epsilon_2)=(\pm,\pm),(+,-)$. The
coalgebraic properties (2.3)-(2.5) may also be expressed concisely: 
\begin{eqnarray}
\bigtriangleup(L^{(\epsilon)})= L^{(\epsilon)}{\dot \otimes} L^{(\epsilon)},
\qquad \varepsilon(L^{(\epsilon)})=1,
\qquad S( L^{(\epsilon)})=(L^{(\epsilon)})^{-1}.
\end{eqnarray}
In the $h\rightarrow 0$ limit, only the `nonstandard' 
deformation parameter $w$ survives and the Hopf algebra ${\cal U}_{h,w}
(h(4))$ reduces to  ${\cal U}_{w}(h(4))$ discussed in [13]. A 
feature of the latter is that the matrices ${\cal R}$ and 
${\tilde {\cal R}}$ become identical as is evident from (2.8) and (2.9) 
in the $h\rightarrow 0$ limit. The above description following the FRT 
procedure then yields the Lax operator $L^{(+)}_{(h\rightarrow 0)}$ and its 
antipode. The matrix generator $L^{(+)}_{(h\rightarrow 0)}$ is given by 
\begin{eqnarray}     
&& L^{(+)}_{(h\rightarrow 0)} =
\pmatrix{1&& 0  & & 0\cr
 && && \cr
 0 && \displaystyle  e^{-wA_{+}} &&\displaystyle   we^{-wA_{+}}N \cr 
 && && \cr
0 && 0&& 1\cr}.
\end{eqnarray}
In the limit $h\rightarrow 0$, the FRT prescription (2.13) and (2.14) 
only describe the Borel subalgebra generated by the elements $A_{+}$, $N$ 
and $E$. 

As the dual Hopf algebra ${\cal F}un_{h,w}({\cal H}(4))$ is also relevant in 
discussing the invariance properties of the RE algebra, we briefly discuss 
it here. Our choice of variables is slightly different from [17]. The function 
algebra ${\cal F}un_{h,w}({\cal H}(4))$ is generated by the variable elements 
of the upper triangular matrix
\begin{eqnarray} 
&&T=\pmatrix{1 & a & b\cr
0 & c & d \cr 
0 & 0 & 1 \cr},
\end{eqnarray}
obeying the relation 
\begin{eqnarray} 
R\;T_{1}\;T_{2}=T_{2} \;T_{1}\;R,
\end{eqnarray}
where $T_{1}=T\otimes 1$ and $T_{2}=1\otimes T$. The algebra reads 
\begin{eqnarray}
&&[a,b]=2h\;a+w\;a^2,\qquad\qquad [c,d]=
w\;c(c-1),\qquad\qquad [b,d]=[a,c ]=0,\nonumber\\
&& [a,d]= w\; a c,\qquad\qquad\qquad 
[b, c]= - w\; a c. 
\end{eqnarray}
Assuming that the diagonal element $c$ is invertible, the coalgebra maps are 
succinctly given by 
\begin{eqnarray}
\bigtriangleup(T)=T{\dot\otimes}T,\qquad\qquad\varepsilon(T)=1,\qquad\qquad
S(T)=T^{-1}. 
\end{eqnarray}
In terms of the generators $a$, $b$, $c$ and $d$, the maps (2.19) read
\begin{eqnarray}  
&& \bigtriangleup (a) = a\otimes c +1\otimes a, \qquad\qquad\;\qquad 
\bigtriangleup(b)=b\otimes 1+1\otimes b+a\otimes d,\nonumber \\
&& \bigtriangleup (c) = c \otimes c,\qquad\qquad\qquad\qquad\qquad
 \bigtriangleup(d)=d \otimes 1+c\otimes d, \nonumber \\
&& \varepsilon (a)=\varepsilon(b)=\varepsilon (d)=0, \qquad\qquad \qquad
\varepsilon (c)=1, \nonumber \\
&& S(a)= - c^{-1} \;a,\qquad S(b)=-b+c^{-1}ad,\qquad
S(c)=c^{-1} , 
\qquad  S(d)= - c^{-1} d.
\end{eqnarray}
The duality of the Hopf algebras ${\cal U}_{h,w}(h(4))$ and 
${\cal F}un_{h,w}({\cal H}(4))$ may now be expressed by the pairing 
\begin{equation}
\langle L^{(\pm)}\otimes 1,\; 1\otimes T\rangle=R_{\pm},
\end{equation}
where $R_{-}=R^{-1}$. 

\sect{The RE Algebra}

Here we study a spectral parameter independent form of RE satisfied by 
the entries of the reflection matrix $K$:
\begin{equation}
RK_1R_{+}K_2=K_2RK_1R_{+},
\end{equation}   
where, as usual, $K_1=K\otimes 1$ and $K_2=1\otimes K$. The entries of 
the matrix $K$ are assumed to be in the upper triangular form   
\begin{eqnarray} 
&&K=\pmatrix{1 & \alpha & \beta\cr
0 &\gamma & \delta \cr 
0 & 0 & 1 \cr},
\end{eqnarray}
A basic property of the RE algebra (3.1) is its covariance under the quantum 
group coaction. The transform $K_{T} =TKT^{-1}$ is also a solution of (3.1) 
if all the elements of the matrices $K$ and $T$ commute 
\begin{eqnarray} 
[K^{i}_{j}, T^{k}_{l}]=0, \qquad\qquad (i,j,k,l)\in (1,2,3). 
\end{eqnarray}
The matrix $T$ satisfies the defining relation (2.17) of the inverse 
scattering 
theory. The RE algebra (3.1) requires its generators $(\alpha,\beta,
\gamma,\delta)$ to have the following commutation properties:
\begin{eqnarray}
&&[\alpha,\beta]=\alpha(2h\;\gamma-w\;\alpha),\qquad\qquad\qquad 
[\alpha,\delta]=(2h\;\gamma-w\;\alpha)(\gamma-1),\nonumber\\
&& [\beta,\delta]=(2h\;\gamma-w\;\alpha)\delta,\qquad\;\;\qquad\qquad 
[\gamma,\;\bullet \;]= 0. 
\end{eqnarray}
The algebra (3.4) has two central elements 
\begin{eqnarray}
C_1=\gamma,\qquad\qquad \qquad C_2=\alpha\delta-\beta\gamma+\beta. 
\end{eqnarray}
The elements $C_1$ is the determinant of the reflection matrix $K$ and is 
assumed to be invertible. This allows us to construct the inverse operator 
$K^{-1}$ as follows 
\begin{eqnarray}
&& K^{-1}=\pmatrix{1 & -\gamma^{-1}\alpha & -\beta+\gamma^{-1}\alpha\delta \cr
0 &  \gamma^{-1} &  -\gamma^{-1}\delta \cr 
0 & 0& 1\cr}.
\end{eqnarray}
Following [30], a representation of the algebra (3.4) may be constructed 
utilizing 
the generators of the quantum algebra ${\cal U}_{h,w}(h(4))$. Using 
the commutation rules (2.1), it may be readily checked that the following 
construction provides a representation of the algebra (3.4):     
\begin{eqnarray}
&& K=S(L^{(-)})L^{(+)} \nonumber \\
&&\phantom{K}=\pmatrix{
1&\displaystyle {2h\over w}e^{-2hE}(1-e^{-wA_{+}})  
& \displaystyle {4h^2\over w}e^{-hE}
(1-e^{-wA_{+}})A_{-}+2he^{-2hE}(1-e^{-wA_{+}})N-2hN \cr
&& \cr
&& \cr
 0 & \displaystyle  e^{-2hE} &
\displaystyle   2he^{-hE}A_{-} -2w e^{-hE}\sinh(hE)N\cr 
&& \cr
&& \cr
0 & 0& 1\cr}\nonumber \\
&& 
\end{eqnarray}
In representation (3.7), the central elements in (3.5) assume the form 
\begin{eqnarray}
C_1=e^{-2hE},\qquad\qquad\qquad C_2=-2he^{-hE}(\sinh(hE)+2hC),
\end{eqnarray}
where $C$ is the Casimir element (2.2) of the ${\cal U}_{h,w}(h(4))$ 
algebra. The representation (3.7) of the $K$-matrix becomes trivial in limit 
$h\rightarrow 0$. This is a consequence of the fact that at $h\rightarrow 0$ 
limit, the universal ${\cal R}$-matrix in (2.8) and its conjugate matrix 
${\tilde {\cal R}}$ in (2.9) become identical. The RE algebra (3.4) in this 
limit, however, remains nontrivial. We will comment on this later.  

As observed in [21] the context of ${\cal U}_{q}(sl(2))$ algebra, the 
construction (3.7) immediately provides a technique to obtain a hierarchy of 
representations of the RE algebra (3.1). Such representations were used in
[21] to 
construct realizations of the braid groups on handlebodies. In the present 
context, successive operations of the coproducts on the representation (3.7)
of the $K$-matrix generate a string of solutions of the RE algebra (3.1):
\begin{eqnarray}
\bigtriangleup(K^{i}_{j}), \qquad (\bigtriangleup\otimes \id)\circ
\bigtriangleup(K^{i}_{j}), \qquad (\bigtriangleup\otimes \id \otimes \id)
\circ (\bigtriangleup\otimes \id)\circ\bigtriangleup(K^{i}_{j}),..., 
\end{eqnarray}
where $(i,j)\in (1,2,3)$. The coproduct used in (3.9) reflects the coproduct 
structure (2.14) of the ${\cal U}_{h,w}(h(4))$ algebra, where
the entries in different sectors are commuting. This contrasts with another 
coproduct scheme [26-28], where the braiding statistics between different 
spaces are to be taken into account. Using (2.14) we may 
successively obtain the representations of the $K$-matrix, whose entries 
now assume values in the tensor product spaces of the 
${\cal U}_{h,w}(h(4))$ algebra. The first few of these representations 
read as follows: 
\begin{eqnarray}  
&& \bigtriangleup({K}^{i}_{j})=(1\otimes S(L^{(-)})^i_l)({K}^l_n\otimes 1)
(1\otimes {L^{(+)}}^n_j),\nonumber \\
&& (\bigtriangleup\otimes \id)\circ
\bigtriangleup({K}^{i}_{j})=(1\otimes 1\otimes S(L^{(-)})^i_p)
(1\otimes S(L^{(-)})^p_l\otimes 1)({K}^l_n\otimes 1\otimes 1 ) \nonumber \\
&&\phantom{ (\bigtriangleup\otimes \id)
\bigtriangleup({K}^{i}_{j})=}(1\otimes {L^{(+)}}^n_m\otimes 1)
(1\otimes 1\otimes {L^{(+)}}^m_j).
\end{eqnarray}

A key property, important from the point of integrability of models, 
requires that two independent solutions of the RE algebra satisfying a 
complementary relation may be combined to construct a new solution of 
the RE algebra. Namely, if $K$ and $K'$ matrices satisfy (3.1), then 
the following combinations 
\begin{eqnarray}  
{\hat K}=KK',\qquad\qquad {\hat {\hat K}}=KK'K^{-1}
\end{eqnarray}
also obey the same RE algebra provided the following complementary relation 
holds:
\begin{eqnarray}  
RK_1R^{-1}K'_2=K'_2RK_1 R^{-1}.   
\end{eqnarray}
The relation (3.12) is covariant under the quantum group coaction. The 
transforms $K_T=TKT^{-1}$ and $K'_T=TK'T^{-1}$ satisfy (3.12) provided 
the elements of $T$ commute with all the elements of $K$ and $K'$. This 
process of building new solutions may obviously be continued. Equation 
(3.12) requires the following commutations relations between the elements 
of $K$ and $K'$:
\begin{eqnarray}
&&\alpha'\alpha=\alpha\alpha', \qquad 
\alpha'\beta=\beta\alpha'-w\;\alpha\alpha',
\qquad \alpha'\gamma=\gamma\alpha',\qquad 
\alpha'\delta=\delta\alpha'-w\;(\gamma-1)\alpha',\nonumber\\
&& \beta'\alpha=\alpha\beta'-2h\;(\gamma-1)\alpha'+w\;\alpha\alpha', \qquad 
\beta'\beta=\beta\beta'-2h\;\delta\alpha'+2hw\;(\gamma-1)\alpha',
 \nonumber \\
&&\beta'\gamma=\gamma\beta',\qquad
 \beta'\delta=\delta\beta'-w\;\delta\alpha'+w^2\;(\gamma-1)\alpha',
\nonumber \\
&& \gamma'\alpha=\alpha\gamma', \qquad
 \gamma' \beta=\beta\gamma'\qquad \gamma'\gamma=\gamma\gamma',\qquad
\gamma'\delta=\delta\gamma',\nonumber \\
&& 
 \delta'\alpha=\alpha\delta'-2h\;(\gamma-1)(\gamma'-1)+w\;\alpha(\gamma'-1),
\qquad 
\delta'\beta=\beta\delta'-2h\;\delta(\gamma'-1)+w\;\alpha\delta',  
\nonumber\\
&&\delta'\gamma=\gamma\delta',  \qquad \delta'\delta=\delta\delta'-w\;\delta(\gamma'-1)+w\;(\gamma-1)\delta'.
\end{eqnarray}
The central elements of $K$ and $K'$ are mutually central in both algebras 
 \begin{eqnarray}
[K^{i}_{j}, C'_{m}]=0,\qquad \qquad [K'^{i}_{j}, C_{m}]=0 \qquad
\hbox{for} \qquad (i,j)\in(1,2,3),\qquad m\in(1,2).
\end{eqnarray}
In (3.13) inhomogeneous terms exist for $(\gamma,\gamma')\not=1$. At this 
point we wish to draw attention to the fact that the RE algebra (3.4) and 
the commutation relations (3.13) between two copies of the generators of the 
algebra (3.4) depend nontrivially on the both the deformation parameters 
$(h,w)$. Relations (3.1) and (3.12) taken together are sometimes mentioned 
as extended RE algebra. The physical meaning of the exchange properties 
(3.13) become 
clearer in the language adopted in [26-28], which we discuss subsequently. 

In an alternate approach [26,27] to the RE algebra, a function type Hopf 
algebra ${\cal B}(R)$, generated by the elements of the $K$-matrix, is 
associated with a $R$-matrix. 
The important distinction is that, in the tensor product structure, the usual 
$\pm$ factor encountered for the twist map of the Hopf superalgebras
is replaced by a map $\psi$ ($\not = \psi^{-1}$) characterizing the braid 
statistics. For $\Gamma_i\in {\cal B}(R)$, where for $i=(1,2,3,4)$, the 
tensor product composition rule reads [26-28]
 \begin{eqnarray}
(\Gamma_1\otimes \Gamma_2)(\Gamma_3\otimes \Gamma_4)=
\Gamma_1\psi(\Gamma_2\otimes \Gamma_3)\Gamma_4,
\end{eqnarray}        
where the braided transposition generates a linear combination of 
tensor product terms. Taking this braiding between two copies of the RE 
algebra into account a new coproduct structure for the $K$-matrices that
is a homomorphism of the RE algebra may be assigned [26-28]. Apart from the 
extra braiding properties, the new Hopf structure is of the same form as for 
the matrix quantum group elements in (2.19):
\begin{eqnarray}
{\tilde \bigtriangleup}(K)=K{\dot\otimes}K,
\qquad\qquad{\tilde \varepsilon}(K)=1,\qquad\qquad
{\tilde S}(K)=K^{-1}. 
\end{eqnarray}  
We repeat here that the braided Hopf structure (3.16) is distinct from the 
Hopf structure of the ${\cal U}_{h,w}(h(4))$ used in (3.9). The latter 
involves no braiding in the tensor product composition rule. 

Expressed in terms of the elements $(\alpha, \beta, \gamma,\delta)$ of the 
$K$-matrix, (3.16) assumes the form
\begin{equation}
\begin{array}{lll}   
{\tilde \bigtriangleup}(\alpha) = \alpha\otimes \gamma +1\otimes \alpha, 
&\qquad {\tilde \varepsilon} (\alpha)=0,
&\qquad {\tilde S}(\alpha)= - \gamma^{-1} \;\alpha, \nonumber \\
{\tilde \bigtriangleup}(\beta)=\beta\otimes 1+1\otimes \beta+\alpha\otimes 
\delta,&\qquad{\tilde \varepsilon}(\beta)=0, &\qquad
{\tilde S}(\beta)=-\beta+\gamma^{-1}\alpha\delta,
\nonumber\\
{\tilde \bigtriangleup}(\gamma) = \gamma \otimes \gamma,
& \qquad {\tilde \varepsilon} (\gamma)=1, &\qquad
{\tilde S}(\gamma)=\gamma^{-1} , \nonumber \\
{\tilde \bigtriangleup}(\delta)=\delta \otimes 1+\gamma\otimes \delta,
&\qquad {\tilde \varepsilon} (\delta)=0 &\qquad
{\tilde S}(\delta)= - \gamma^{-1} \delta.
\end{array}
\end{equation}
The braiding rules between two copies of the RE algebra (3.4) may be directly 
read from (3.13). The primed quantities are entries in the second sector of 
the tensor product space and quantities like $\Gamma'_1\Gamma_2$, where
$(\Gamma_1,\Gamma_2)\in (\alpha,\beta,\gamma,\delta)$, are interpreted as a 
way of writing the transpose $\psi(\Gamma_1\otimes \Gamma_2)$. In these 
notations (3.13) reads
\begin{eqnarray}
&& \psi(\alpha\otimes\alpha)=\alpha\otimes \alpha, \qquad 
\psi(\alpha\otimes\beta)=\beta\otimes \alpha
-w\;\alpha\otimes \alpha, \nonumber \\
&&\psi(\alpha\otimes\gamma)=\gamma\otimes \alpha, \qquad 
\psi(\alpha\otimes\delta)=\delta\otimes \alpha-w\;(\gamma-1)\otimes \alpha
, \nonumber \\
&& \psi(\beta\otimes\alpha)=\alpha\otimes \beta-2h\;(\gamma-1)\otimes \alpha+
w\;\alpha\otimes \alpha, \nonumber \\
&& \psi(\beta\otimes\beta)=\beta\otimes \beta-2h\;\delta\otimes \alpha
+2hw\; \;(\gamma-1)\otimes \alpha,\nonumber \\
&& \psi(\beta\otimes\gamma)=\gamma\otimes\beta,
\qquad \psi(\beta\otimes\delta)=\delta\otimes \beta-w\;\delta\otimes \alpha
+w^2\;(\gamma-1)\otimes \alpha, \nonumber \\
&& \psi(\gamma\otimes\alpha)=\alpha\otimes \gamma,
\qquad \psi(\gamma\otimes\beta)=\beta\otimes \gamma,\qquad
\psi(\gamma\otimes\gamma)=\gamma\otimes \gamma,
\qquad \psi(\gamma\otimes\delta)=\delta\otimes \gamma, \nonumber \\
&& \psi(\delta\otimes\alpha)=\alpha\otimes \delta-2h\;(\gamma-1)\otimes 
(\gamma-1)+w\;\alpha\otimes (\gamma-1), \nonumber \\
&& \psi(\delta\otimes\beta)=\beta\otimes \delta-2h\;\delta\otimes (\gamma-1)
+w\; \alpha\otimes \delta,\nonumber \\
&& \psi(\delta\otimes\gamma)=\gamma\otimes\delta,
\qquad \psi(\delta\otimes\delta)=\delta\otimes \delta-w\;\delta\otimes
(\gamma-1)+w\;(\gamma-1)\otimes \delta, 
\end{eqnarray}
The braided exchange properties (3.18) in the $w\rightarrow 0$ limit was 
earlier discussed in [28] the context of the braided group corresponding 
to the `standard' deformation of the Heisenberg algebra ${\cal U}_h(h(4))$.
The central elements $(C_1,C_2)$ in (3.5) are truly bosonic in nature in 
the sense
\begin{eqnarray}
\psi(C_m\otimes \Gamma)=\Gamma\otimes C_m\qquad\hbox{for} \qquad \Gamma
\in {\cal B}(R)\qquad\hbox{and}\qquad m\in(1,2).
\end{eqnarray}  
The elements $C_2$ is, however, not grouplike
\begin{eqnarray}
{\tilde \bigtriangleup}(C_2)\not= C_2\otimes C_2.
\end{eqnarray} 
Taking into account, the braiding properties (3.18), it may be checked that 
the maps $({\tilde \bigtriangleup},{\tilde S},{\tilde \varepsilon})$ 
preserve the algebra (3.4). In this construction the antipodes of the 
composite operators are defined by the rule [26-28]
 \begin{eqnarray}
{\tilde S}(\Gamma_1\Gamma_2)=m\circ\psi({\tilde S}(\Gamma_1)\otimes 
{\tilde S}(\Gamma_2))
\end{eqnarray} 
where $m$ is the multiplication map. As an example, we demonstrate the property
\begin{eqnarray}
{\tilde S}([\beta,\delta])={\tilde S}(2h\;\gamma\delta-w\;\alpha\delta).
\end{eqnarray} 
Using (3.18) and (3.21), we obtain: 
\begin{eqnarray}
&& \hbox{l.h.s.}=m\circ(\gamma^{-1}\psi(\beta\otimes \delta)
-\psi(\delta\otimes \beta)
\gamma^{-1}+\gamma^{-1}\psi(\delta\otimes \alpha\delta)\gamma^{-1}-
\gamma^{-1}\psi(\alpha\delta\otimes \delta)\gamma^{-1}) \nonumber \\
&& \phantom{l.h.s.}=-2h\;\gamma^{-2}\delta-w\;\gamma^{-2}\alpha\delta-
2wh\;\gamma^{-1}+2wh\nonumber \\
&& \hbox{r.h.s}=-2h\;\gamma^{-2}\delta-w\;m\circ(\gamma^{-1}
\psi(\alpha\otimes\delta)\gamma^{-1}) =\hbox{l.h.s.} \nonumber 
\end{eqnarray}

We now examine the RE algebra (3.4) in two distinct limits. In the limit  
$w\rightarrow 0$, the doubly deformed algebra ${\cal U}_{h,w}(h(4))$ 
reduces to the algebra ${\cal U}_{h}(h(4))$ with the
standard deformation parameter $h$. The RE algebra (3.4) now assumes the 
form of a centrally extended Lie algebra:
\begin{eqnarray}
&&[\alpha,\beta]=2h\;\alpha\gamma,\qquad\qquad [\alpha,\delta]=
2h\;\gamma(\gamma-1),\qquad\qquad
[\beta,\delta]=2h\;\gamma\delta,
\end{eqnarray}
where the central element $\gamma$ equals to a number in an irreducible 
representation. When $\gamma\not =0$, the remaining generators $(\alpha,\beta,
\gamma)$ may be rescaled as follows:      
\begin{eqnarray}
J={\beta\over 2h\gamma}, \qquad\qquad P_+=\delta,\qquad\qquad P_-=\alpha.
\end{eqnarray} 
The generators $(J,P_{\pm})$ then satisfy the algebra:
\begin{eqnarray}
[J,P_{\pm}]=\pm P_{\pm}, \qquad\qquad [P_+,P_-]=\xi, 
\end{eqnarray} 
where the central charge $\xi=2h\gamma(1-\gamma)$ vanishes when $\gamma=1$. 
The invariant $C_2$
of the RE algebra in (3.5) becomes the Casimir element of the algebra (3.25) 
and reads 
\begin{eqnarray}
C_2=P_-P_+ +\xi J.
\end{eqnarray} 
The algebra (3.25) is, in fact, isomorphic to that of the magnetic translation 
vectors in two dimension, where the central charge term appears because 
of a magnetic flux perpendicular to the plane. In the $\gamma\not=1$ case, 
the RE algebra (3.23) therefore indicates an effective magnetic flux, which 
may be of relevance in discussing statistical properties. The braided 
commutations relations between two copies of the algebra (3.23) may be 
obtained from (3.13) at $w\rightarrow 0$ limit
\begin{eqnarray}
&&\alpha'\alpha=\alpha\alpha', \qquad 
\alpha'\beta=\beta\alpha',\qquad 
\alpha'\delta=\delta\alpha',\nonumber\\
&& \beta'\alpha=\alpha\beta'-2h(\gamma-1)\alpha', \qquad 
\beta'\beta=\beta\beta'-2h\delta\alpha',\qquad \beta'\delta=\delta\beta',
 \nonumber \\
&&\delta'\alpha=\alpha\delta'-2h(\gamma-1)(\gamma'-1), \qquad 
\delta'\beta=\beta\delta'-2h\delta(\gamma'-1), 
\qquad \delta'\delta=\delta\delta'.
\end{eqnarray}
The central terms $\gamma$ and $\gamma'$ commute with all others elements 
of the extended RE algebra.

In the `nonstandard' $h\rightarrow 0$ limit, only the dimensional deformation
parameter $w$ is present. The Hopf algebra ${\cal U}_{h,w}(h(4))$ 
reduces to ${\cal U}_{w}(h(4))$ considered in [13] and the RE algebra 
(3.4) now reads  
\begin{eqnarray}
[\alpha,\beta]=-w\;\alpha^2,\qquad\qquad [\alpha,\delta]=
-w\;(\gamma-1)\alpha, \qquad\qquad
 [\beta,\delta]= -w\; \alpha \delta,
\end{eqnarray}
where the central charge element $\gamma$ is assumed to have a numerical 
value. Assuming $\alpha$ to be the invertible, the relations (3.28) 
may be expressed as a linear algebra as follows. For an invertible $\alpha$, 
we choose a set of a new variables 
\begin{eqnarray}
X=\alpha,\qquad\qquad Y=\alpha^{-1}\beta,\qquad\qquad Z=\delta
\end{eqnarray}
to reexpress (3.28) as a Lie algebra 
\begin{eqnarray}
[X,Y]=-wX,\qquad [X,Z]=-w(\gamma-1)X,\qquad 
[Y,Z]=w(\gamma-1)Y-wZ
\end{eqnarray}
The Casimir element $C_2$ now reads
\begin{eqnarray}
C_2=X(Z-(\gamma-1)Y).
\end{eqnarray}
We distinguish between two cases:

\smallskip

\smallskip

{\it (1).} When $\gamma=1$, we redefine ${\hat J}=\displaystyle {Y\over w}$, 
${\hat P}_{+}=X$, ${\hat P}_{-}=Z$. The algebra (3.30) is now isomorphic to 
the inhomogeous algebra $iso(1,1)\simeq t_2 \odot so(1,1)$:
\begin{eqnarray}
[{\hat J},{\hat P}_{\pm}]=\pm {\hat P}_{\pm}, \qquad\qquad [{\hat P}_+,
{\hat P}_-]=0
\end{eqnarray} 
with its Casimir invariant $C_2={\hat P}_{+}{\hat P}_{-}$. 

\smallskip

\smallskip

{\it (2).} When $\gamma\not=1$, we use the variables  
\begin{eqnarray}
{\hat {\hat J}}={1\over 2w}(Y+{Z\over \gamma-1}), 
\qquad\qquad 
{\hat {\hat P}}_+={X\over w^2\;(\gamma-1)}
\qquad\qquad 
{\hat {\hat P}}_-={1\over w}(Y-{Z\over \gamma-1})
\end{eqnarray} 
The algebra (3.30) is again isomorphic to $iso(1,1)$: 
\begin{eqnarray}
[{\hat {\hat J}}, {\hat {\hat P}}_{\pm}]=\pm 
{\hat {\hat P}}_{\pm}, \qquad\qquad [{\hat {\hat P}}_+,{\hat {\hat P}}_-]=0
\end{eqnarray} 
with the corresponding invariant operator 
\begin{eqnarray}
C_2=-w^3(\gamma-1)^2 {\hat {\hat P}}_{+}{\hat {\hat P}}_{-}.
\end{eqnarray} 
Thus, in the `non-standard' ($h\rightarrow 0$) limit, the RE algebra (3.28) 
may be identified with a Lie algebra  $iso(1,1)$ provided $\alpha$ is 
invertible. As a real form, the algebra $iso(1,1)$ is isomorphic to 
(1+1)-dimensional Poincar\'e algebra. A deformed $iso(1,1)$ algebra was 
previously obtained [31] at a
contraction limit of the Jordanian deformation of the $sl(2)$ algebra.
The unitary representations of the algebra $iso(1,1)$ are infinite dimensional.
The braiding relations for the exchanges between two copies of the RE algebras 
(3.28) may be read from (3.13) at the `nonstandard' $h\rightarrow 0$ limit:  
\begin{eqnarray}
&&\alpha'\alpha=\alpha\alpha', \qquad 
\alpha'\beta=\beta\alpha'-w\;\alpha\alpha',\qquad 
\alpha'\delta=\delta\alpha'-w\;(\gamma-1)\alpha',\nonumber\\
&& \beta'\alpha=\alpha\beta'+w\;\alpha\alpha', \qquad 
\beta'\beta=\beta\beta',
 \qquad 
\beta'\delta=\delta\beta'-w\;\delta\alpha'+w^2\;(\gamma-1)\alpha', \\
&&  \delta'\alpha=\alpha\delta'+w\;\alpha(\gamma'-1), \qquad
\delta'\beta=\beta\delta'+w\;\alpha\delta', \nonumber\\
&& \delta'\delta=\delta\delta'-w\;\delta(\gamma'-1)+w\;(\gamma-1)
\delta'.\nonumber 
\end{eqnarray}

The RE algebra (3.4) corresponding to the doubly-deformed $(h\not=0, w\not=0)$
$R$-matrix (2.12) is essentially nonlinear in character and cannot be reduced 
to a linear algebra.

\sect{Conclusion}

In summary, we considered the RE algebra for a finite dimensional $R$-matrix 
for the doubly-deformed algebra ${\cal U}_{h,w}(h(4))$. A representation of the reflection matrix $K$ was constructed using the matrix generators $L^{(\pm)}$
 of the ${\cal U}_{h,w}(h(4))$ algebra. The coproduct rules of the 
${\cal U}_{h,w}(h(4))$ algebra then yields a hierarchy of solutions of the 
$K$-matrix. The complementary condition necessary for combining two distinct 
solutions of the RE algebra yields the braiding relations between the 
generators of these two algebras. This may be thought as a generalization 
of Bose-Fermi statistics to braiding statistics, which then may be used 
to provide a new braided coalgebraic structure to a Hopf algebra generated 
by the elements of matrix $K$. The RE algebra and the braiding relations 
between two copies of the RE algebra depend on both the `standard' and the 
`non-standard' deformation parameters $h$ and $w$ respectively. For 
$w\rightarrow 0$ the RE algebra is found to be isomorphic to the 
algebra of magnetic translation operators on a two-dimensional plane. For the 
`non-standard' $(h\rightarrow 0)$ limit, the RE algebra becomes isomorphic 
to the inhomogeneous algebra $iso(1,1)$, which is the Poincar\'e algebra 
in $(1+1)$ dimension. When both deformation parameters are present $(h\not=0,
w\not=0)$, the RE algebra (3.4) is non-linear in nature.

\vskip 1cm

\noindent {\bf Acknowledgments:} 

\smallskip 

We thank Amitabha Chakrabarti and Ramaswamy Jagannathan for interesting 
discussions. One of us (RC) wants to thank Amitabha Chakrabarti for a 
kind invitation. He is also grateful to the members of the CPTH group 
for their kind hospitality.

\bibliographystyle{amsplain}

\end{document}